\documentclass[a4paper,11pt]{article}
\usepackage{xypic,amssymb}
\title{ Joins for (Augmented) Simplicial Sets}
\author{P. J. Ehlers  and T. Porter,\\School of Mathematics,\\University of Wales Bangor,\\Dean Street, \\ Bangor, Gwynedd, LL57 1UT,\\ Wales, U.K.}

\newtheorem{example}{Example}[section]
\newtheorem{theorem}[example]{Theorem}
\newtheorem{proposition}[example]{Proposition}

\newtheorem{lemma}[example]{Lemma}
\newcommand{\qed}{\hfill $\square $ }
\begin{document}
\maketitle
 \begin{abstract}
We introduce a notion of join for (augmented) simplicial sets generalising
the classical join of geometric simplicial complexes.  The definition comes
naturally from the ordinal sum on the base simplicial category $\Delta$.

{\noindent\bf 1991  Math. Subj. Class.:} 18G30
\end{abstract}
\section{Introduction}

The theory of joins of (geometric) simplicial complexes as given by Brown, \cite{Top:Brown}, or Spanier, \cite{:Span}, reveals the join operation to be a basic geometric construction.  It is used in the development of several areas of geometric topology (cf. Hudson, \cite{Hudson}) whilst also being applied to the basic properties of polyhedra relating to homology.

The theories of geometric and abstract simplicial complexes run in a largely parallel way and when describing the theory, expositions often choose which aspect -- abstract  combinatorial or geometric -- to emphasise at each step.  Historically in algebraic topology geometric simplicial complexes, as tools, were largely replaced by CW complexes whilst the combinatorial abstract complex became part of simplicial set theory.  In the process, joins were negelected and there does not seem to be a well known definition of the join of two simplicial sets.

Within the setting of simplicial set theory, the ordinal sum plays a strange r\^ole. This operation takes two ordinals and concatenates them, so $[m]or[n] = [m+n+1]$, where \linebreak$[m] = \{ 0 < 1 < \cdots < m \}$, so it is fundamental for the combinatorics of ordinals.  In the literature on simplicial set theory it seems rarely to be mentioned, yet it is sometimes there but hidden, for instance, in the $\overline{W}$-construction for simplicial groups (see May, \cite{Algtop:May}) or simplicial groupoids (see Dwyer and Kan, \cite{DandK}).  

In this context it occurs through the use of the Artin-Mazur codiagonal, \cite{AandM}, which assigns to a bisimplicial set or group, a much smaller model of the homotopy type than does the diagonal.  (The diagonal is intuitively easier to use and tends to be ``wheeled out'' whenever passage from bisimplicial objects to simplicial objects is needed; however, it may not always be the most efficient tool to use.)  This codiagonal is linked with the total DEC functor (Illusie \cite{Cot:Ill}, Duskin \cite{:Dus}, Porter, \cite{:Porter}, Bullejos et al \cite{Bull}), which can be given explicitly in terms of the ordinal sum.  

In this brief note, it is shown that the ordinal sum leads naturally to a ``join''  operation on {\em augmented}  simplicial sets, and the relation of this join to the geometric join is studied.

\section{Definitions}\label{defn}

It will be assumed that the reader is conversant in general with basic
simplicial set theory, in particular, the definition of the
 singular complex of a topological space, and the geometric realisation of a singular complex.  On the subject of notation, note that the simplicial set which is called the $n$-simplex, $\triangle [n]$\glossary{$\triangle [n]$}, is the representable functor, $\Delta (-,[n])$. The simplicial set $\triangle [n]$ will be referred to as the {\it standard} $n$-{\it simplex}.

The category of finite ordinals will be denoted  $\Delta $:  the
ordinal $\{0 < 1 < \cdots < n\}$ will be denoted $[n]$ with the empty
set being denoted by $[-1]$.\\

\noindent {\bf Definition \ref{defn} (i)}\\
Let $f_i:[p_i] \rightarrow [q_i]$ for $i = 0, 1$. Define the
``ordinal sum'' functor, \\ $or:\Delta ^{2} \longrightarrow
\Delta$, as follows:-\[ or([p_0],[p_1]) = [p_0 + p_1 + 1]  \]\[
or(f_0,f_1) = \left\{\begin{array}{lcl} f_0(k) & \mbox{ if } & 0 \leq
k \leq p_0 \\ f_1(k-p_0-1)+q_0+1 & \mbox{ if } & p_0+1 \leq k
\\\end{array} \right. \] \\
Note that $[-1]$ is a two sided identity for the operation on
objects.\\

\noindent {\bf Definition \ref{defn} (ii)} \\
An augmented simplicial set is a simplicial
set, $X$, together with an augmentation, that is, a set $X_{-1}$ and
a morphism $q_X:X_0\longrightarrow X_{-1}$, where $q_Xd_0 = q_Xd_1$.

There is an obvious forgetful functor from the category, $ASS$, of augmented simplicial sets
to that $SS$, of simplicial sets, cf. Duskin \cite{:Dus}.\\

\noindent {\bf Definition \ref{defn} (iii)} \\
The {\it canonical augmentation} of a simplicial set has $X_{-1} =
\pi_0(X)$ and $q_X$ the coequaliser of 
$$X_1 \begin{array}{l}
\stackrel{d_0}{\longrightarrow} \\ \stackrel{d_1}{\longrightarrow} 
\end{array} X_0 .$$
This augmentation is left adjoint to the forgetful functor.\\

\noindent {\bf Definition \ref{defn} (iv)} \\
The {\it trivial augmentation} of a simplicial set has $X_{-1} =
\ast $, the one point set, and $q_X$
the unique (trivial) morphism  $X_0 \longrightarrow \ast$.
This augmentation is right adjoint to the forgetful functor.\\

\noindent {\bf Definition \ref{defn} (v)} \\
The {\it geometric realisation} defined on augmented simplicial sets is the
composition of the forgetful functor to simplicial sets and the usual 
geometric realisation functor to topological spaces.

(This is the only reasonable definition of a geometric realisation on
augmented simplicial sets, as the codomain of the augmentation is, in
some sense, the image of the empty set.)

\noindent {\bf Definition \ref{defn} (vi)} \\
The {\it singular complex functor} from topological spaces to
augmented simplicial sets is the composition of the normal singular complex
functor, which is right adjoint to the geometric realisation functor,
with the trivial augmentation functor, right adjoint to the forgetful functor.

It is automatic that the two functors so defined are adjoint.

\section{Combinatorial Join}
\label{monoid}

The following is our proposed definition for a join of augmented simplicial sets.\\

\noindent {\bf Definition \ref{monoid} (i)}\\
Let the {\it join} of two augmented simplicial sets $X$ and $Y$ be denoted $X \odot Y $.  The set of $n$-simplices,  $(X \odot Y)_n$,  is:-
\[ \bigsqcup_{i=-1}^n X_{n-1-i} \times Y_i \]
the face maps are given by:- 
\[d_i^n(x,y) = \left\{ \begin{array}{cl} (d_i^px,y) & \mbox{ \ if \ } 0
\leq i \leq p \\(x,d_{i-p-1}^{n-p-1}y) & \mbox{ \ if \ } p < i \leq n
\end{array} \right. \]
where $(x,y) \in X_p \times Y_{n-p-1}$, and $d_0^0$ is the
augmentation (of $X$ or $Y$);\\
lastly, the degeneracies are:-
\[s_i^{n-1}(x,y) = \left\{ \begin{array}{cl} (s_i^px,y) & \mbox{ \ if \ }
0 \leq i \leq p \\(x,s_{i-p-1}^{n-p-2}y) & \mbox{ \ if \ } p < i \leq n-1
\end{array} \right. \]
where $(x,y) \in X_p \times Y_{n-p-2}$.

There is also a coend definition for $\odot$:\[ X \odot Y
\cong\int^{p,q} (X_p \times Y_q) \cdot \triangle([p]or[q]) \]

\noindent {\bf Remark}\\
It is trivial to prove that  $\triangle [n] \odot \triangle [m]
\cong \triangle [n+m+1]$.  It is also true that $\odot $ is an
associative operation, but the simplest proof requires a number of
constructions and results associated with the join which are not of
immediate interest here.

\section{Topological Join}
\label{top}
The following  definition is a generalisation of the concept of join for two
suitable subspaces of a vector space.
The {\it topological join} thus defined is discussed in some detail in chapter 5,
section 7 of \cite{Top:Brown}. Results proved there will be used here
without proof: the notation for this section is largely taken from
there.  We work within the category of compactly generated spaces.\\

\noindent {\bf Definition \ref{top} (i)} \\
 Consider two topological spaces
${\cal U}$ and ${\cal V}$, and construct a set  of 4-tuples
$(r,u,s,v)$, where $u \in {\cal U},\; v \in {\cal V}, \; r,s \in
[0,1]$ and $r + s = 1$: in
the case that $r = 0$, the $u$ will be ignored, and in the case that
$s = 0$, the $v$ will be ignored. This set will be suggestively called
${\cal U} \ast {\cal V}$.

There are obvious projections from this set of 4-tuples: \\
$p_{\cal U}:{\cal U} \ast {\cal V} \rightarrow {\cal U}$,
\rule{5pt}{0pt} $p_{\cal V}:{\cal U} \ast {\cal V} \rightarrow {\cal
V}$, \rule{5pt}{0pt} $p_r:{\cal U} \ast {\cal V} \rightarrow (0,1] $
and $p_s:{\cal U}
\ast {\cal V} \rightarrow (0,1]$ \\
which are termed the {\it coordinate functions} of ${\cal U} \ast
{\cal V}$. The first two are obviously defined, the last two take a
point $(r,u,s,v) \in {\cal U} \ast {\cal V} $ to $r$ and $s$
respectively.

The {\it topological join} of ${\cal U}$  and  ${\cal V}$ is
defined as the set ${\cal U} \ast {\cal V}$ together with the initial
topology with respect to the {\it coordinate functions}. Thus a
function with codomain ${\cal U} \ast {\cal V}$ is a continuous
function if and only if its composite with each of the coordinate
functions is continuous. 

\medskip 

To compare the combinatorial and topological join operations, we will need more precision on the construction of the geometric realisation.
There are a number of different constructive definitions of geometric
realisation. The process is essentially the following: 

(i) take one copy of $\triangle ^n$ for each non-degenerate $n$-simplex of
$X$;\\
and then 

(ii) glue them all together using the face and degeneracy maps of
the simplicial set $X$ (see \cite{Cat:Mac}). \\ Explicitly we have: 

Let $X$ be a simplicial set. Define $RX$ by:\\
\[ RX = \sqcup_{n \in {\mathbb N}}\sqcup_{x \in X_n} \triangle ^n_x \]
Define an equivalence relation on $RX$ as generated by the
following relation:\\ 
writing $({\bf p},x)$ for $(p_0,\cdots ,p_m) \in \triangle ^m_x$ 
and $({\bf q},y)$ for $ (q_0,\cdots ,q_n) \in \triangle ^n_y$ then
$({\bf p},x) \sim ({\bf q},y)$ if either \\
\rule{20mm}{0mm} $d_ix = y$ and $\delta _i(q_0,\cdots ,q_n) = 
(p_0,\cdots , p_m)$ or \\
\rule{20mm}{0mm} $s_ix = y$ and $\sigma _i(q_0,\cdots ,q_n) = 
(p_0,\cdots, p_m)$,\\
where the $\delta_i$ and $\sigma_i$ are the  continuous maps given by face inclusion and folding in the usual way.
Then $ |X| \cong {RX}/{\sim} $ where $ {RX}/{\sim } $ has 
the identification topology.

 \begin{proposition}
\label{join} 
\[\triangle ^p \ast \triangle ^q \cong \triangle ^{p+q+1} \]
\end{proposition}
{\bf Proof} \\
Consider the vector space ${\mathbb R}^{p+q+1}$ and the two compact convex subsets:
\[ X = \{ \, (x_0,x_1,\cdots, x_p,0,\cdots ,0)
\, | \, \sum_{i=0}^p x_i = 1 \, \} \]
\[ Y = \{ \, (0,\cdots,0,y_0,y_1,\cdots ,y_q) \, | \, \sum_{j=0}^q
y_j = 1 \, \} \]
First note that $X \cong \triangle ^p$ and  $Y \cong \triangle ^q$.
Furthermore, it is clear that no two lines in the set $U =  \{ \,
r{\bf x} + (1-r){\bf y} \, | \, 0 \leq r \leq 1, \, {\bf x} \in X, \,
{\bf y} \in Y \, \} $ intersect except at endpoints. Thus $X \ast Y =
U$. However, $U$ is the subset of ${\mathbb R} ^{p+q+1} $ given by
\[ \{\, (rx_0,\cdots ,rx_p, (1-r)y_0,\cdots ,(1-r)y_q) \, | \,
\sum_{i=0}^p rx_i + \sum_{j=0}^q (1-r)y_j = 1 \, \}. \]
That is, $U$ is the affine $(p+q+1)$-simplex. Therefore $\triangle ^p \ast
\triangle ^q \cong \triangle ^{p+q+1}$. \qed \\

When we form $\triangle [p] \odot \triangle [q]$, we obtain, on varying $p$ and $q$, a bicosimplicial object in $SS$.  (In general if $\cal C$ is a category, a cosimplicial object in $\cal C$ is a functor from $\Delta$ to $\cal C$, whilst a bicosimplicial object is a functor from $\Delta \times \Delta$ to $\cal C$.) Similarly $\triangle^p \ast \triangle^q$ is a bicosimplicial space.

 \begin{lemma}
There is a natural isomorphism
\[ |\triangle [p]| \ast |\triangle [q]| \cong |\triangle [p] \odot
\triangle [q]| \]
of bisimplicial spaces.
\label{join2}
\end{lemma}
{\bf Proof}\\
Recall $|\triangle [m]| := \triangle ^m $. Since 
$\triangle [p] \odot \triangle [q] \cong \triangle ([p]or[q]) = 
\triangle [p+q+1]$, the isomorphism exists for each pair $(p,q)$.
Now $\{\triangle ^n\}_{n \in  {\mathbb N}}$ has an obvious cosimplicial 
structure, and  the isomorphism is easily seen to be an isomorphism of bicosimplicial 
  spaces. \qed\\

 \begin{theorem} \rule{0pt}{12pt} \\
Let $X$ and $Y$ be trivially augmented simplicial sets. Then
\[ |X \odot Y| \cong |X| \ast |Y| \]
\label{join3}
\end{theorem}
{\bf Proof}\\ 
(The following is a direct geometric proof: we will comment later on the categorical aspects.)

  Recall that 
\[ |X| \ast |Y| \, :=  \,  \left\{ \begin{array}{ll}  r[(p_0,\cdots,p_m)_x]
\rule{24pt}{0pt} s.t. & \sum_{i=0}^m p_i = 1 \mbox{ , } \sum_{i=0}^n
q_i = 1 \mbox{ , } \\ \rule{12pt}{0pt} + s[(q_0, \cdots q_n)_y] \,  &
x \in X_m \mbox{ , } y \in Y_n \mbox{ , } r + s = 1  
\\ & p_i,q_i,r,s \geq 0 \\ & \mbox{and }
[-] \mbox{ denotes equivalence class} \end{array} \right\} \]

It should also be noted that if $r = 0$, the point from $|X|$ is
ignored and similarly, if $s = 0,$ the point from $|Y|$ is ignored.\\
Define a map $ f:|X| \ast |Y| \longrightarrow |X \odot Y| $ as
follows:
\[ f(r[(p_0,\cdots ,p_m)_x] + s[(q_0,\cdots ,q_n)_y]) \mapsto
[(rp_0,\cdots ,rp_m,sq_0,\cdots ,sq_n)_{x,y}] \]

The function $f$ is  well defined since if $r = 0$, the point $x$ is ignored, similarly if $s = 0$. This means that for any $y$, it
must be true that  $(0,\cdots 0,q_0,\cdots ,q_n)_{(x,y)} \sim (0,\cdots
0,q_0,\cdots q_n)_{(x',y)}$ for all $x,x' \in X$. This will be true exactly when the augmentation of both $X$ and $Y$ are trivial as was required. A moment's thought then will show that the function $f$
respects the relation and so is well defined.  Continuity is also
trivial to check. The obvious inverse function is also continuous under the
definition of the topology on $|X| \ast |Y|$. Thus the two spaces are homeomorphic. \qed\\

\noindent{\bf Remarks}\\

 (i)  It may seem slightly contrived that the condition ``trivially augmented'' should be needed, however consider the following example:

Let  $X := \triangle [0] \sqcup \triangle [0]$ together with the {\it canonical} augmentation, and consider $X \odot X$.  The result is the disjoint union of four unit intervals -- that is, $\triangle [1] \sqcup \triangle [1] \sqcup \triangle [1] \sqcup \triangle [1] $:  Ideally, the result should be homotopically equivalent to a $1$-sphere.

(ii)  The Theorem above is in fact a simple consequence of a categorical argument which shows a different aspect of the necessity for having a trivial augmentation.

The simplicial complex functor to augmented simplicial sets needs to specify an augmentation, and for the functor to be right adjoint to the geometric realisation functor, the augmentation must be the trivial one (since the trivial augmentation is right adjoint to the forgetful functor from augmented simplicial sets to simplicial sets).  Thus the condition `trivially augmented' merely requires that the augmented simplicial sets are related to the geometric realisation functor upon which the theorem depends. The result is now seen to depend just on left adjoints interacting nicely with the coends in the geometric realisation and  join functors.
\section{Simplicial Spheres.}
Recall (from \cite{Top:Brown}) that
\[ {\bf S}^p \ast {\bf S}^q \cong {\bf S}^{p+q+1} . \]
This essentially says that the $n$-sphere in the category of topological spaces is the join of $n+1$ copies of the $0$-sphere.

 There are several simplicial models for the n-sphere.  For instance, Gabriel and Zisman, \cite{GandZ}, p.26, define the simplicial circle, $\Omega$, to be the coequaliser of the pair of morphisms
$$\diagram
 \triangle [0]\rto<1ex>^{\delta_0}\rto<-1ex>_{\delta_1} & \triangle[1],
\enddiagram$$
and the suspension of a pointed simplicial set $X$ to be $\Omega \wedge X$.  This gives an n-sphere as being $\bigwedge^n\Omega$, obtained from the n-cube $\triangle[1]^n$ by collapsing the `boundary' of the cube to a point.  Other authors form a simplicial sphere by collapsing the boundary $\partial \triangle[n]$ of the n-simplex to a point.

The join operation suggests another form.
Consider the simplicial set formed as the disjoint union of two copies of $\triangle [0]$ and augmented trivially.  This will be denoted by ${\bf S}^0$  and will be referred to as the simplicial $0$-sphere.
Then  ${\bf S}^0 \odot {\bf S}^0$ has four non-degenerate $1$-simplices connected to each other in a ``diamond'' as below:-\\
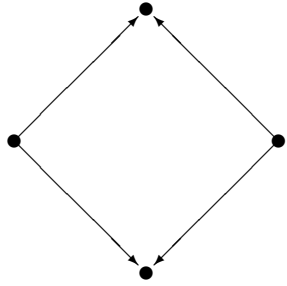
\begin{picture}(300,115)
 
\put(150,55){\vector(1,1){47}}
\put(150,55){\vector(1,-1){47}}
\put(250,55){\vector(-1,1){47}}
\put(250,55){\vector(-1,-1){47}}
\put(150,55){\circle*{5}}
\put(200,105){\circle*{5}}
\put(200,5){\circle*{5}}
\put(250,55){\circle*{5}}
\end{picture}

Define the simplicial $n$-sphere, ${ S}^n \in ob{\mathit ASS}$, as follows:-
\[S^n : = \underbrace{ {S}^0 \odot \cdots \odot { S}^0}_{n+1}  \]
It is clear from the definition of combinatorial join and of the simplicial $0$-sphere that the simplicial $n$-sphere is a triangulation of the topological $n$-sphere.  In fact, theorem~\ref{join3} gives explicitly that 
\[ |{ S}^n| \cong {\bf S}^n. \]  
Moreover this model clearly satisfies 
\[ S^p \odot  S^q \cong S^{p+q+1}\]
 unlike the other models.  Thus if we write $\Sigma^n = \triangle[n] / \partial\triangle[n]$ then $\Sigma^p \odot \Sigma^q$ has one non-degenerate simplex in each of the dimensions 1, p + 1, q + 1, and p + q + 1, and two non-degenerate simplices in dimension 0 and so `looks' totally unlike $\Sigma^{p+q+1}$.

The combinatorial join forms part of a closed monoidal structure on the category of augmented simplicial sets, $ASS$. (The `internal hom' is given by 
\[ [X,Y]_n = ASS(X,Dec^{n+1}Y),\]
where $Dec$ is the d\'ecalage functor (see Duskin, \cite{:Dus} ).)  It is therefore possible to define augmented analogues of the loopspace construction that are compatible with the join.

\end{document}